\documentclass[12pt,a4paper]{article}
\usepackage{amsmath,amsthm,amsfonts,amssymb}
\def\ni{\noindent}

\theoremstyle{definition}

\def\l{\lambda}

\begin{document}

\begin{center}
 \Large{\bf Structure of the tensor product of two  simple modules of quantum $GL_2$}\\
\vspace{1cm}
 \small{ M. S. Datt\\
{\it  School of Mathematics \& Statistics\\
 University of Hyderabad, P.O Central University, Hyderabad, India-500046.\\
 e-mail: msdatt@uohyd.ac.in}}
 
\end{center}

 Abstract: In this article, we consider the tensor product of two simple modules of quanum $GL_2$ over a field of characteristic $p\neq 0$. We show that it can be expressed as a direct sum of indecomposable twisted tilting modules. This problem has been studied by Henke and Doty [1] for $SL_2$ and also later on for $SL_3$ by S R Doty, Chris Bowman and Stuart Martin, ([7, 8]).\\

\ni {\it Key words: Simple Modules, Tilting modules, Twisted tilting modules.}

\section{Preliminaries}
 In this section, we describe terminology and notation. We fix a field $k$ of characteristic $p\neq 0$. Let $q$ be a primitive $\ell^{th}$ root of unity in $k$. By  a quantum $k$-group $G$,  we mean a Hopf algebra $k[G]$ over $k$.  We are concern herewith the quantum $GL_n$ as introduced by R. Dipper and S. Donkin in [3]. Let $A_q(n)$ denote the $k$- algebra generated by $x_{ij}$, $1\leq i,j\leq n$ subject to the relations (i) $x_{ir}x_{is}=x_{is}x_{ir}$ for all $1\leq i,r,s,\leq n$ (ii) $x_{is}x_{jr}=q^{-1}x_{jr}x_{is}$ for all $1\leq i<j\leq n$ and $1\leq r\leq s\leq n$ and finally $x_{js}x_{ir}=x_{ir}x_{js}+(q-1)x_{is}x_{jr}$ for all $1\leq i<j\leq n$. The algebra $A_q(n)$ is a biialgebra with the comodule structure maps $x_{ij}\mapsto \sum_{r=1}^n x_{ir}\otimes x_{rj}$ (comultiplication) and counit $\epsilon(x_{ij})=\delta_{ij}$ (the kronecker delta). Let $S_n$ be the group of permutations on $\{1,2,3\cdots ,n\}$ and $d_q=\sum _{\sigma\in S_n}sgn(\sigma)x_{i,i\sigma}x_{2,2\sigma}\cdots x_{n,n\sigma}$, the quantum determinant. The localisation $A_q(n)_{d_q}$ of $A_q(n)$  at $d_q$ (as $\{1, d_q, d_q^2\cdots \}$ is a Ore set) is a Hopf algebra. By the quantum $GL_n(k)$ we mean the Hopf algebra $A_q(n)_{d_q}$. From now on we denote the quantum  $GL_n$ over $k$ by $G$ and its corresponding Hopf algebra by $k[G]$. We also call the Hopf algebra $k[G]$ as the coordinate ring of $G$. We denote the coordinate functions of $k[G]$ by $c_{ij}$, $1\leq i,j\leq n$. By a quantum subgroup $H$ of $G$ we mean the Hopf algebra $k[G]/I_H$ for some Hopf ideal $I_H$ of $k[G]$. We are concern herewith two important quantum subgroups of $G$, namely a Borel subgroup of $G$ and torus contained in $B$.  Let $I_B$ be the Hopf ideal of $k[G]$ generated by $\{c_{ij}~ | ~1\leq i<j\leq n\}$ and the  corresponding quantum subgroup of $G$ is called Borel subgroup of $G$. The quantum subgroup $T$ with the corresponding Hopf ideal $I_T$ of $k[G]$ generated by $\{c_{ij}~|~i\neq j\}$ is called torus. By definition, we have $T\subset B\subset G$. Let ${\rm mod}(H)$ denote the category of finite dimesnional $G$-modules, for $H=G,B,T$. We need to discuss some representation theory of quantum $GL_n(k)$ as we need in the sequel. 

We let $X(T)={\mathbb{Z}}^n$ and ${\mathbb{Z}}X(T)$ be the group algebra with basis $\{e(\lambda)~|~\lambda\in X(T)\}$ and multiplication is given by the rule $e(\l)e(\mu)=e(\l+\mu)$. By [7], for each $\lambda=(\lambda_1, \lambda _2\cdots , \lambda_n)\in X(T)$ we have the $T$-module $k_{\lambda}$, with the comodule structure map $a\mapsto a\otimes c_{11}^{\lambda_1}c_{22}^{\lambda_2}\cdots c_{nn}^{\lambda_n}$. Also the set $\{k_{\lambda}~|~\lambda\in X(T)\}$ is a complete set of mutually non-isomorphic simple $T$-modules. Any finite dimensional $T$- module is completely reducible. For  $V\in {\rm mod}(T)$, let $V(\lambda)$ be the sum of all submodules of $V$ isomorphic to $k_{\lambda}$. We denote the formal character $\chi(V)$ of $V$ by $\sum_{\lambda\in X(T)}V(\lambda)e(\lambda)$.  We say that an element $\l=(\l_1,\l_2,\cdots ,\l_n)\in X(T)$ is a dominant weight if $\l_1\geq \l_2\geq \cdots \l_n$. We denote the set of all dominant weights by $X^+(T)$.
	
We have a homomorphism $\phi: B\longrightarrow T$, whose comorphism $\widetilde{\phi}: k[T]\longrightarrow k[B]$ takes $c_{ii}$ to $c_{ii}+$ the defining ideal of $B$. Hence any $T$- module can be treated as a $B$-module via $\phi$. For each $\l\in X(T)$, we have the induced module $\nabla(\l)={\rm Ind}_B^Gk_{\l}$. The module $\nabla(\l)\neq 0$ if and only if $\l\in X^+(T)$.
For any $\l\in X^+(T)$,  $\nabla(\l)$ is called costandard modules. The socle $L(\l)$  of $\nabla(\l)$ is simple, for any $\l\in X^+(T)$. Moreover, the set $\{L(\l)~|~\l\in X^+(T)\}$ is a complete set mutually non-isomorphic simple $G$-modules.  For $\l\in X^+(T)$, let $\Delta(\l)=\nabla(-w_0(\l))^*$, where $w_0$ is the longest element of $S_n$ and $\nabla(-w_0(\l))^*$ is the dual of $\nabla(-w_0(\l))$. The modules $\Delta(\l)$ are called standard/Weyl modules.
 For $V\in {\rm mod}(G)$, we call a filtration $0\subset V_0\subset V_1\subset \cdots \subset V_r\subset\cdots  $ of submodules of $V$ is a good filtration, if $V=\cup_{i=1}^{\infty}V_i$, 
and for each $i\geq 1$, we have  $V_i/V_{i-1}$ is either 0 or isomorphic to $\nabla(\lambda)$,  for some $\lambda \in X^+(T)$.  Similarly, if a module $V\in {\rm mod}(G)$ has a filtration $(0)=V_0\subset V_1\subset \cdots \subset V_r\subset \cdots$ of submodules of $V$ such that the quotient $V_i/V_{i+1}$ is either $(0)$ or isomorphic to $\Delta(\l)$ for some $\l\in X^+(T)$, then we say that $V$ has a standard filtration.  A module $V\in {\rm mod}(G)$ is said to be a tilting module if $V$ has a good filtration as well as  a Weyl filtration. 

By a theorem of C. M.  Ringel [10], for each  $\lambda\in X^+(T)$ we have a unique(upto isomorphism) indecomposable tilting module $T(\lambda)$ of highest weight $\lambda$. The set $\{T(\lambda)~|~\lambda\in X^+(T)\}$ is a full set of mutually non-isomorphic indecomposable tilting modules. These tilting modules are called   partial tilting modules. And also any  tilting module is a direct sum of partial tilting modules $T(\lambda)$, $\l\in X^+(T)$. To discuss representation theory of quantum $GL_n(k)$, we also need to recall some representation theory of general linear group $GL_n$ over $k$.

We write $\overline{G}$ for the algebraic group $GL_n(k)$. Let $y_{ij}$ be the coordinate functions for $\overline{G}$. Let $\overline{B}$ be a Borel subgroup of $\overline{G}$  containing a torus  $\overline{T}$.  Let $X^{+}(\overline{T})$ be the set of all dominants weights of $\overline{G}$. Let $\overline{F}: \overline{G}\longrightarrow \overline{G}$ be the Frobenius morphism(ordinary). 
By [3, 1.3.2 Corollary] we have the Frobenius morphism $F: G\longrightarrow \overline{G}$ whose comorphism takes $y_{ij}$ to $c_{ij}^{\ell}$. Hence any $\overline{G}$-module can be treated  as $G$ module via $F$. In a paper,  by Ann Henke and S R Doty [1], it is shown that the tensor product of two simple $SL_2$- modules can be decomposed as a direct sum of twisted tensor product of certain tilting modules. In this paper, we show that the tensor product of two simple quantum $GL_2$-modules is a direct sum of certain indecomposable modules. And each indecomposable summand is the tensor product of  a certain tilting module with the Frobenius twist($F$) of twisted tilting module of $\overline{G}$. Our methods are tilting modules, Steinberg's tensor product theorem for quantum groups, Clebsch-Gordan formula. 

\section{The Quantum $GL_2.$}
In this section, we shall specialize to quantum $GL_2(k)$.
Let $G$ be the quantum $GL_2(k)$. 
Let $X_1=\{(a,b)\in X^+(T)~|~0\leq a-b\leq \ell -1, a, b\geq 0\}$ and $\pi=\{(a,b)\in X^+(T)~|~a,b\geq 0, ~0\leq a-b\leq 2(\ell-1)\}$. The set $\pi$ is a saturated subset of $X^+(T)$ in the sense that for $\l\in \pi$ and $\mu\in X^+(T)$, if $\mu\leq \l$ then $\mu\in \pi$. Therefore  the category of finite dimensional $G$-modules $V$ whose composition factors  belonging to $\{L(\l)~|~\l\in \pi\}$ is equivalent to the category of finite dimensional modules for the generalized $q$-schur algebra(see [9]).  We have the following Lemma. \\

\ni {\bf Lemma 2.1 } {\it (a) If $\l\in X_1$, then $L(\l)=\Delta(\l)=\nabla(\l)=T(\l)$.\\
	(b) Let $\l\in \pi\setminus X_1$ and $\l=\ell+r$, $0\leq r\leq \ell-2$. Then the composition factors of $T(\l)$ are isomorphic to $L(\l)$ or $L(\l-(r+1)(\epsilon_1+\epsilon_2)$ and the character of $T(\l)$ is given by $\chi(\l)+\chi(\l-(r+1)(\epsilon_1-\epsilon_2)$.}\\

\ni {\it Proof:} (a) This follows by  [7] or [11].

\ni (b) We follow the arguments given in [7]. Let $\l\in \pi\setminus X_1$. By [8, 3.2], we know that $\nabla(\l)$ is uniserial. 
Let $\l=(a,b)\in \pi,~ a-b=\ell +r$,  $0\leq r\leq \ell -2$ and  $\mu=\l-(r+1)(\epsilon _1-\epsilon _2)$.   Then by [7], the injective hull $I(\l)$ of $L(\l)$ in  mod$(\pi)$ is $\nabla(\l)$. Now by [3], we have $((I(\l):\nabla(\tau))=[\nabla(\tau), L(\l)]$. As $\l\in \pi\setminus X_1$, $(I(\mu): L(\l'))=1$ if $\l'=\l-(r+1)(\epsilon _1-\epsilon_2)$ or $\l' =\l$. By [2], the block in mod$(\pi)$ containing $\l$ is $\{\l, \lambda-(r+1)(\epsilon _1-\epsilon _2)\}$. Therefore $T(\l)$ has $\nabla$-filtration with sections $\nabla(\l), \nabla(\l-(r+1)(\epsilon _1+\epsilon_2))$.

Let St$=(\ell -1+b,b)$ be the Steiberg weight and $\mu =(a-b-(\ell -1), 0)$.  Let $\l= {\rm St}+\mu$. By (a), St and $L(\mu)$ are tilting. Hence by [7] the module  $\nabla({\rm St})\otimes \nabla(\mu)$ is tilting  and has  highest weight $\l$. Therefore the module  $T(\lambda)$ will be the block component of $\nabla({\rm St})\otimes \nabla(\mu)$ of highest weight $\l$ [4 Prop(1.2)]. By [7], we note that two tilting modules are isomorphic if and only if they have same character.
By Clebsch-Gordan formula, we have $\chi({\rm St})\chi(\mu)= \chi(\l)+\chi(\l-(\epsilon _1-\epsilon _2)\cdots +\chi (\l- (r+1)(\epsilon _1-\epsilon _2))$. As any  tilting module is a direct sum partial tilting modules, by putting altogether, we see that  $\chi(\l)+\chi(\l-(r+1)(\epsilon_1-\epsilon _2))$ is the character of $T(\l)$. \hfill{$\Box$}\\

Let $\mathbb{T}=\{T(\l)~|~\l\in \pi\}$ and we call them as special tilting modules of quantum $GL_2$.  By   [5, 2.2] any tilting module of $SL_2$ can be expressed as a direct sum of twisted tensor product of special tilting modules.   This remark still holds for the general linear group  $GL_2$ as well as quantum $GL_2$ also(see [7]). In the case of quantum $GL_2(k)$, we  call these  modules as $(\overline{F})^F$-twisted tilting modules of $G$. In Theorem 3.2, we show that the tensor product $L\otimes L'$ of two simple module $L,L'$ of quantum $GL_2$ can expressed as a direct sum of $(\overline{F})^F$ twisted tilting of $G$. All of these summands are indecomposable. In general, a direct summand of this expression  need not tilting. But $L,L'$ are tilting  if and only if every indecomposable summand is tilting. 
To prove this, it is enough consider the  dominant weights of the form $(a,0)$. We can see this as follows: let $\l=(a,b)\in X_1$ and $d_q=L(\epsilon_1+\epsilon _2)$. We can write $\l=(a-b,0)+b(\epsilon _1+\epsilon_2)$. Then we have $L(\l)=d_q\otimes L(a-b,0)=d_qL(a-b,0)$.  

Let $\l=(a, 0), \mu=(c, 0)\in X_1$.  Then  $a+b<\ell$ or $a+b=\ell +r,$ $r\leq \ell -2$. By Clebsch-Gordan formula, we have $\chi(\l)\chi(\mu)=\chi(a+b,0)+\chi(a+b-1,1)+\cdots +\chi(a,b)$ or $\chi(\l)\chi(\mu)=\sum _{i=0}^{r+1}\chi(a+b)-i(\epsilon_1+\epsilon _2)$, if $a+b= \ell+r$. We let $W(\l, \mu)=\{(a+b,0), (a+b-1,1),\cdots (a,b)-(r+1)(\epsilon_1+\epsilon_2)\}$ and $W_S(\l,\mu)=\{(c,d)\in W(\l ,\mu)~|~c-d\geq \ell\}$. We need the following lemma.\\

\ni {\bf Lemma 2.2(2):} {\it Let $(a,0), (b,0)\in X_1$ and $L=L(a,0)$ and $L'=L(b,0)$ be two simple modules of highest weights $(a,0)$ and $(b,0)$ respectively.  Then $L\otimes L'$ is a tilting module and it is isomorphic to the direct sum of tilting modules $T(\lambda)$, $\lambda$ varies over $W(L, L')\setminus W_S(L, L')$.}\\

\ni {\it Proof:} We know that the tensor product of two tilting modules is a tilting module, therefore  $L\otimes L'$ is a tilting module.  Now by the Lemma 2.1, if $W(L,L')\subset X_1$, then $L\otimes L'$ is isomorphic to the direct sum $\oplus _{\gamma \in W(L,L')}T(\gamma)$. Suppose  $W(L, L')\not\subset X_1$. If $(c,d)\in W_S(L, L')$, then as $L\otimes L'$ is tilting, $\chi(c,d)+\chi((c,d)-(c-d+1)(\epsilon _1+\epsilon _2))$ is the character of $T(c,d)$, by the Lemma 2.1. Hence $L\otimes L'$ is isomorphic to the direct sum $\oplus _{\gamma\in W(L,L')\setminus W_S(L,L')}T(\gamma)$.\hfill{$\Box$} \\

\ni {\bf Lemma 2.3(3):} {\it Let $(a,0), (b,0)\in X_1$ and let $L=L(a,0)$ and $L'=(b,0)$, where $(a,0),(b,0)\in X_1$. Then the tilting moddule $L\otimes L'$ is  indecomposable if and only if $a=0$ or $((a,b)=(\ell-1, 1)$.}\\

\ni {\it Proof:} If $a=0$ then $L\otimes L$ is isomorphic to $L'$ and it is indecomposable by the Lemma 2.1. In case of $(a,b)=(\ell-1, 1)$, the comoposition factors of $L(\ell-1,0)\otimes L(1,0)$ are $L(\ell ,0)$ and $L(\ell-(\epsilon _1+\epsilon _2))$. The modules  $L(\ell-1,0)$ and $L(1,0)$ are tilting and hence $L(\ell-1,0)\otimes L(1,0)$ is also tilting. Now by Lemma 1.2, the  character of $L(\ell-1,0)\otimes L(1,0)$ is $\chi(\ell, 0)+\chi(\ell-(\epsilon _1+\epsilon _2))$ and this is same as the character of $T(\ell, 0)$ and hence it is indecmposable tilting, as required.\\

Conversely, suppose $L\otimes L'$ is tilting. Suppose $(a,b)\neq (\ell-1, 1)$.  First we consider the case  $a=\ell-1$ and $b>1$. In this case we have  $a+b=\ell +r$, $1\leq r \leq \ell-2$. By Clebsch-Gordan formula, $\sum _{i=1}^{r+1}\chi (a+b-i(\epsilon_1+\epsilon _2)$ is the character of $L\otimes L'$ . From this we can conclude that  the indecomposbale tilting modules $T(\ell+1, 0)$ and $T(\ell-1, 1)$ appear as direct summands in the decompostion of $L\otimes L'$ as direct sum of indeomposable tilting modules, which is contradiction. Similarly, we can conclude that $L(a,0)\otimes L(b,0)$ cannot be indecomposbale if $b=1$ and $a< \ell-1$ \hfill{$\Box$}\\

\ni {Example} (i) For $\ell =5$, we have $\ell-1=4$. Let $a=4$ and $b=2$. Then  we have $L(4, 0)\otimes L(2,0)=T(6, 0)\oplus T(5,1)$.\\
(ii) Suppose $\ell=5$. Let $a=3$ and $b=1$. Then $L(3,0)\otimes L(1,0)=T(4,0)\oplus T(3,1)$.\\

\ni For any $\l\in X_1$ and $\mu \in X^+$, by the Steinberg's tensor product theorem [See [7]], we hve $L(\l)\otimes L(\mu)^{F})\simeq L(\l+\ell\mu)$,  where $F$ is the Frobenius (quantum) morphism. Given any $a\in {\mathbb{N}}^+$, we can write  $a=\tau+\ell(\sum_{i=0}^ra_ip^i)$, $0\leq \tau\leq \ell-1$ and $0\leq a_i\leq p-1$, uniquely. We call this as the $(\ell, p)$ expnasion of $a$. Now for $\l=(a,0)$ and $\mu=(\ell\sum_{i=0}^ra_ip^r, 0)$, we have $L(a,0)\simeq L(\tau, 0)\otimes L(\sum _{i=0}^{r}a_i p ^i , 0)$.  Since  $0\leq a_i\leq p-1$,  we have $\overline{L}(a_i,0)=\overline {T}(a_i, 0)$ for all $i$ and similarly $L(\tau ,0)=T(\tau ,0)$ as $0\leq \tau\leq \ell-1$. Hence we have the following Lemma.\\

\ni {\bf Corollary 2.4} {\it Suppose $a\in {\mathbb{Z}}^+$ and $a=\tau +\ell (\sum_{i=0}^ra_ip^ i)$ be the $(\ell, p)$ expansion of $a$. Then $L(a, 0)=L(\tau)\otimes (\otimes _{i=0}^r \overline\displaystyle{{\overline{L(a_i)}}^{\overline{F}^ i})}^F$.}\hfill{$\Box$}\\

\section{ Structure of the tensor product of Simple Modules}

In this section, we show that the tensor product of two simple $G$-modules can be expressed as a finite direct sum of $F$-twist of $\overline{F}$-twisted tensor product of special tilting modules. As we discussed above, we prove this for simple $G$-modules of the form $L(a,0)$, $a\in {\mathbb{Z}}^+$. Let $\overline{X_1}=\{(a,b)~|~0\leq a-b\leq p-1\}$ and $\overline{\pi}=\{(a,b)\in X^+({\overline{T}})~|~0\leq a-b\leq 2(p-1)\}$.  First we prove the following Lemma:
\\ 

\ni {\bf Lemma 3.1} {\it For $M\in {\rm mod}(G)$ and $N\in {\rm mod}(\overline{G})$, we have $M\otimes N^F\simeq N^F\otimes M$.}\\

\ni {\it Proof:} Let $\{x_1,x_2,\cdots , x_m\}$ and $\{y_1,y_2,\cdots ,y_n\}$ be bases for $M$ and $N$ respectively. We shall show that $M\otimes N^F\simeq N^F\otimes M$. Let $\psi_M: M\longrightarrow M\otimes k[G]$, $\psi(x_i)=\sum _{r} x_r\otimes f_{ri}$ (resp. $\psi_N (y_j)=\sum_{s} y_s\otimes g_{sj}$) be the structure map for M (resp. for $N$). Let $\varphi : M\otimes N^F\longrightarrow N^F\otimes M$ be the  map defined as $\varphi(x\otimes y)=y\otimes x$.  

Let $\psi$ be the structure map for $M\otimes N$, i.e., $\psi(x\otimes y)=\sum _{r,s} x_r\otimes y_s\otimes f_{ri}g_{sj}$. Now consider $(\varphi\otimes 1)(\psi(x_i\otimes y_j))=(\varphi\otimes 1)(\sum x_r\otimes y_s\otimes f_{ri}g_{sj})=\sum_{r,s}y_s\otimes x_i\otimes f_{ri}g_{sj}$. The module $N^F$ is a $G$-module with the structure map given by $\psi_{N^F}(y_j)=\sum_{s}y_s\otimes F(g_{sj})$. By [3], we have $F(f)$ is in the center, for any $f\in k[G]$. 

Let  $z_j$ be $y_j$ when we viewed as an element of $N^F$. Hence $\psi(\varphi\otimes 1)(x_i\otimes z_j)=\sum _{r,s}z_{s}\otimes x_{r}\otimes F(g_{sj})f_{ri}=\sum_{r,s}z_{s}\otimes x_r\otimes F(g_{sj})f_{ri}$. On the otherhand, we have,
$(\varphi\otimes 1)\psi (x_i\otimes z_j)=(\varphi\otimes 1)(\sum _{r,s}x_r\otimes z_s\otimes f_{ri}F(g_{sj})=\sum_{r,s}z_s\otimes x_r\otimes F(g_{sj})f_{ri}$. Hence $M\otimes N^F\simeq N^F\otimes M$.\hfill{$\Box$}\\

\ni {\bf Theorem 3.2 } {\it Let $a,b\in {\mathbb{Z}}^+$. Then the tensor product $L(a,0)\otimes L(b,0)$ can be expressed is a finite direct sum of indecomposable modules of the form} $$T(\lambda _{-1})\otimes (\overline{T}(\l_0)\otimes \overline{T}(\l_1)^{\overline{F}}\otimes \cdots \otimes \overline{T}(\l_r)^{\overline{F}^r})^F$$

\ni where $\l_{-1}\in X_1$ and $\l_{i}\in \overline{\pi}$, for $0\leq i\leq r$\\

\ni {\it Proof:} Let $a=\tau _a+\ell (\sum _{i=0}^ma_ip^i)$ and $b=\tau _b+\ell (\sum _{j=0}^n b_jp^j)$ be  $(\ell, p)$ expansion of $a$ and $b$ respectively. Without loss of generality we can assume that $m=n$.  By the the Steinberg's tensor product theorem, we have $L(a,0)=L(\tau_a, 0)\otimes (\otimes_{i=0}^m\overline{L}(a_i,0)^{\overline{F}^i} )^F$ and $L(b,0)=L(\tau _b,0)\otimes (\otimes _{i=0}^n \overline{L}(b_i,0)^{\overline{F ^i}})^F$. By the Lemma 3.1, we have $L(a,0)\otimes L(b,0)=(L(\tau _a, 0)\otimes L(\tau _b,0))\otimes (\otimes _{i=0}^n ((L(a_i,0)\otimes L(b_i, 0))^{\overline{F}^i})^F.$ 

Asuume that $\tau_a\geq \tau_b$. By Clebsch-Gordan formula, the character of  $L(\tau_a, 0)\otimes L(\tau _b, 0)$ is given by $\chi(\tau_a+\tau_b,0)+\chi(\tau_a+\tau_b-1,1)+\cdots +\chi(\tau_a,\tau _b)$. Now as $(\tau_a, 0), (\tau_b,0)\in X_1$, by Lemma 2.2, we have $L(\tau_a,0)\otimes  L(\tau_b,0)$ is a direct sum of tilting modules $T(\l)$, where $\l\in W(L(\tau_a,0),L(\tau_b,0))\setminus  W_S(L(\tau_a,0),L(\tau_b,0))$. Similarly, as the weights $(a_i,0),(b_i,0)\in \overline{X}_1$,  we can express  $\overline{L}(a_i,0)\otimes \overline{L}(b_i,0)$ as a direct sum of tilting modules $\overline{T}(\mu)$, where  $\mu$ varies over $\overline{W}(L(a_i,0),L(b_i,0))$, for all $0\leq i\leq n$. 

Let $L=L(\tau _a,0)$, $L'=L(\tau _b, 0)$ and $\overline{L}_i=\overline{L}(a_i,0)$, $\overline{L}_i'=\overline{L}(b_i,0)$,and $I=W(L,L')\setminus W_S(L,L')$ and $I_i=W(L_i,L_i')\setminus W_S(L.L')$, for $0\leq i\leq n$. Then we have $L\otimes L'=\oplus_{\l\in I}T(\l)$ and $\overline{L}_i\otimes \overline{L}'_i=\oplus _{\mu\in I_i}\overline{T}(\mu)$, for $0\leq i\leq n$. Thus we have,  
$$L(a,0)\otimes L(b,0)=(\oplus_{\l\in I}T(\l))\otimes (\otimes _{i=0}^n(\oplus _{\mu\in I_i}\overline{T}(\mu))^{\overline{F}^i})^F$$

By interchanging tensor prodcut with direct sums, we can express $L(a,0)\otimes L(b,0)$ can be expressed as a direct sum of modules of the form $$T(\lambda _{-1})\otimes (\overline{T}(\l_0)\otimes \overline{T}(\l_1)^{\overline{F}}\otimes \cdots \otimes \overline{T}(\l_r)^{\overline{F}^r})^F$$

By [$q$-Schur], we see that each term in the sum is indecomposable.\\

\ni Let $\underline{\l}=(\l_{-1}, \l_{0},\l_{1},\cdots \l_{n})$, where $\l_{-1}\in \pi$ and $\l_{i}\in \overline{\pi}$. Let $\l_i=(a_{i1},a_{i2})$, for $-1\leq i\leq r$. Let $M(\underline{\l})=T(\l_{-1})\otimes (\otimes_{i=0}^n\overline{T}(\l_i)^{\overline{F}^r})^F$.  \hfill{$\Box$}\\




\ni {\bf Acknowledgements:} {\it It gives me great pleasure to thank my teacher S. Donkin for suggesting the problem and valuable guidance throughout the work. }


\begin{thebibliography}{99}
\bibitem{} Ann Henke and S. R. Doty {\it Decomposition of tenor product of modular irreducibles for $SL_2$,}  Quarterly Jurnal of Mathematics, 56, (2005), 189-207. 
\bibitem{} A. G. Cox {\it On Some applications of infinitesimal methods to quantum groups and related algebra}, Ph.D. thesis, University of London, 1997.
\bibitem{} R. Dipper and S. Donkin {\it Quantum $GL_n$}, Proc. London Mathematical society (3), 63 (1991), 165-211.
\bibitem{} S. Donkin {\it The Blocks of a semisimple algebraic Group}, Journal of Algebra, Vol. 67, No. 1 November 1980.
\bibitem{} S. Donkin {\it On tilting modules for algebraic groups}, Math. Zeit. 212 (1993), 39-60.
\bibitem{} S. Donkin, {\it standard homological properties for quantum $GL_n$}, J. Algebra, 181, (1996), 235-266,
\bibitem{} S. Donkin {\it The $q$-Scur Algebra}, Cmbridge University Press, London Mathematical Society Lecture Note Series. 253, 1998.
\bibitem{} S R Doty, Chris Bowman and Stuart Martin, {\it Decomposition of tenor product of modular irreducibles for $SL_3$,  Journal of Algebra 9, (2011), 177-219.  }
 \bibitem{} S R Doty, Chris Bowman and Stuart Martin, {\it Decomposition of tenor product of modular irreducibles for $SL_3$, the $p\geq 5$ case}, Int. Electron, Journal of Algebra, 17, (2018), 105-138.\\
\bibitem{} C. M. Ringel, {\it The category of modules with good filtrations over a quasi-hereditary algebra has almost split sequences}, Math. Zeit. 208 (1991), 209-225.
\bibitem{} L. Thamas, {\it The subcomodule structure of quantum symmetric powers}, Bull. Australian Math. Soc. 50 (1994), 29-39.

\end{thebibliography}
\end{document}